# ON THE TESTABILITY OF THE CAR ASSUMPTION

By Eric A. Cator[1]

*Delft University of Technology*

In recent years a popular nonparametric model for coarsened data is an assumption on the coarsening mechanism called coarsening at random (CAR). It has been conjectured in several papers that this assumption cannot be tested by the data, that is, the assumption does not restrict the possible distributions of the data. In this paper we will show that this conjecture is not always true; an example will be current status data. We will also give conditions when the conjecture is true, and in doing so, we will introduce a generalized version of the CAR assumption. As an illustration, we retrieve the well-known result that the CAR assumption cannot be tested in the case of right-censored data.

**1. Introduction.** When dealing with coarsened data, the coarsening may be due to some random effect. A condition was proposed in Heitjan and Rubin (1991) on this random effect, called "coarsened at random," or CAR. In their setup the random variable of interest, which in this paper we will call $Y$, takes values in a finite set $\mathcal{Y}$. However, instead of observing $Y$ directly, we observe a nonempty random set $X \subset \mathcal{Y}$ such that with probability 1, $Y \in X$. They then define the CAR assumption as an assumption on the possible or allowed conditional distributions of $X$ given $Y = y$ [CAR is a modelling assumption, so a class of distributions for $(Y, X)$ is considered]:

for all $A \subset \mathcal{Y}$     $\mathbb{P}(X = A | Y = y)$ is constant in $y \in A$.

They showed that in this setting, the CAR assumption ensured that the randomness of the coarsening could be ignored when making inference on the parameter of interest, namely, the distribution of $Y$. Many papers have since appeared generalizing this idea, especially to general sample spaces. We refer to Jacobsen and Keiding (1995) and Gill, van der Laan and Robins (1997) for a general introduction. Our goal is mainly to discuss the testability of

Received July 2002; revised December 2003.
[1]Supported in part by Netherlands Organization for Scientific Research (NWO).

the CAR assumption, that is, does the CAR assumption restrict the possible distributions of the data?

We will start by giving a general model for coarsened data which is very close to the one given in Jacobsen and Keiding (1995), but without the measurability issues in that paper. We repeat that it is not our main goal to extend the notion of CAR to general sample spaces; therefore, we will not give an extensive comparison with definitions given in the aforementioned papers. We would just like to mention that in practical situations all definitions will lead to more or less the same concept. Furthermore, our notation will mostly be similar to that in Pollard (2002), with one notable exception: if $\mu$ is a measure on a space $\mathcal{Z}$, and $\pi$ is a measurable map from $\mathcal{Z}$ to $\mathcal{Y}$, then we denote the image measure on $\mathcal{Y}$ of $\mu$ under $\pi$ as $\pi(\mu)$.

Let $\mathcal{Y}$ be the space of the variable of interest $Y$ (e.g., the time of onset of a certain disease). The stochastic variable $Y$ is distributed according to a probability measure $Q$. Let $\mathcal{Z}$ be a "hidden" space from which we can retrieve $Y$ and the data. To be more precise, the stochastic variable $Z \in \mathcal{Z}$ is distributed according to a probability measure $\mu$ and there exists a measurable map $\pi:\mathcal{Z} \to \mathcal{Y}$ such that $Y = \pi(Z)$. Furthermore, there exists a measurable map $\psi:\mathcal{Z} \to \mathcal{X}$, where $\mathcal{X}$ is the data space, such that $X = \psi(Z)$ is the observed data. In short,

$$\begin{array}{c} (\mathcal{Z},\mu) \xrightarrow{\pi} (\mathcal{Y},Q). \\ \psi \downarrow \\ (\mathcal{X},P) \end{array}$$

The measure $\mu$, together with the mappings $\pi$ and $\psi$, contains all the information about how the variable of interest $Y$ is *coarsened* into the data $X$. This definition of coarsened data is more general than the one used by, for example, Gill, van der Laan and Robins (1997), where the data must consist of sets. However, it is also much easier to find counterexamples to the conjecture mentioned in the abstract, to which we will come shortly.

First, to make things a bit more tangible, let us see how current status data fits into our framework: let $Y$ be the time of onset of a certain disease, let $C$ be the time of visiting a doctor, generally called the censoring time, and define the data $X$ as

$$X = (C, \mathbb{1}_{\{Y \leq C\}}).$$

Then $Z = (Y,C)$ (so $\mathcal{Z} = [0,\infty[ \times [0,\infty[)$, $\pi(Y,C) = Y$ and $\psi(Y,C) = (C, \mathbb{1}_{\{Y \leq C\}})$.

In Heitjan and Rubin (1991), Gill, van der Laan and Robins (1997), Nielsen (2000) and several others, coarsened data consists of sets $B$, elements of some $\sigma$-algebra $\mathcal{B}$ on $\mathcal{Y}$, such that $Y \in B$. Defining $\mathcal{Z} = \mathcal{Y} \times \mathcal{B}$, we see that this approach also fits into ours if we have proper conditions on $\mu$: we allow all $\mu$ such that $Y \in B$ almost surely. Of course, we could also say that our data consists of the set $\pi(\psi^{-1}\{x\}) \subset \mathcal{Y}$; it is, however, possible (see Example 2.1) that knowing $x$ provides more information. In any case, we find that our results are more clearly stated in our definition of coarsened data.

Before we can state the CAR assumption, we need some more notation. We will restrict ourselves in this paper to dominated models, so we choose a fixed and known probability measure $\mu_0$ on $\mathcal{Z}$. In Gill, van der Laan and Robins (1997) the CAR assumption is also defined for the nondominated case [CAR(ABS)], but we will get back to this later. Define

$$Q_0 = \pi(\mu_0) \quad \text{and} \quad P_0 = \psi(\mu_0).$$

Now we wish to condition on the map $\pi$ (or, equivalently, on $Y$). If $\mathcal{Z}$ and $\mathcal{Y}$ are, for example, Polish spaces, this can always be done via a Markov-kernel: we define the conditional distribution of $Z$ under $\mu_0$ given $Y = y$, denoted by $\mu_0(dz|y)$, such that for each bounded measurable function $k$ on $\mathcal{Z}$ we have

$$\int_{\mathcal{Z}} k(z)\mu_0(dz) = \int_{\mathcal{Y}} \left( \int_{\mathcal{Z}} k(z)\mu_0(dz|y) \right) \pi(\mu_0)(dy).$$

This is called a disintegration. Of course, we also have that

$$\mu_0(\{z : \pi(z) \neq y\}|y) = 0$$

for $\pi(\mu_0)$-almost all $y$.

DEFINITION 1.1 (The CAR assumption). In the notation given above, the CAR assumption states that $\mu \ll \mu_0$ is a possible (or admitted) distribution of $Z$ if and only if

$$\mu(dz) = g \circ \psi(z) \cdot h \circ \pi(z) \mu_0(dz),$$

where $h$ is an arbitrary density with respect to $Q_0$ and $g$ is a positive measurable function on $\mathcal{X}$ such that

(1.1) $$\int_{\mathcal{Z}} g \circ \psi(z) \mu_0(dz|y) = 1 \qquad \text{for } Q_0\text{-almost all } y,$$

which is equivalent with

$$E_{\mu_0}(g(X)|Y) = 1.$$

This implies that $h(y)$ is the (marginal) density of $Y$ with respect to $Q_0$ and that the conditional distribution of $Z$, given $Y = y$, is given by

$$\mu(dz|y) = g \circ \psi(z) \mu_0(dz|y).$$

This loosely means that we assume that given $Y$, the *unknown* part by which the coarsening mechanism chooses $Z$ (note that $\mu_0$ is known!) may only be a function of the data. Note that under CAR, we can choose an arbitrary density $h \in L^1(Q_0)$, but the measurable function $g$ must be positive and satisfy (1.1) [in particular, $g \in L^1(P_0)$]. This restriction on $g$ does not depend on $h$, however, which gives the set of all possible distributions of $Z$ under CAR a product structure.

It might not be entirely clear why one would want to make such an assumption, but the popularity of the CAR assumption can largely be explained by the following proposition. First, we define a linear map

(1.2) $$S : L^1(Q_0) \to L^1(P_0) : S(h)(x) = E_{\mu_0}(h(Y)|X=x)$$

(remember that if $Z \sim \mu_0$, then $X \sim P_0$ and $Y \sim Q_0$).

PROPOSITION 1.2. *Let $\mu$ be a distribution of $Z$ that satisfies the CAR assumption. This means that there exists $g \in L^1(P_0)_+$ such that $\mu(dz|y) = g \circ \psi(z)\mu_0(dz|y)$. Let $h$ be the marginal density of $Y$ with respect to $Q_0$ [so $\pi(\mu)(dy) = h(y)Q_0(dy)$]. Then the marginal distribution of $X$ is given by*

$$\psi(\mu)(dx) = g(x)S(h)(x)P_0(dx).$$

This shows that the likelihood of the data factorizes into a relevant factor [remember that $h(Y)$ as a function of $h$ is the likelihood based on the underlying data $Y$, the variable of interest, and note that $S$ is known] and a nuisance factor $g$. Since we can choose any $g$ that satisfies (1.1) and then choose an arbitrary density $h$ independent of the chosen $g$, the overall parameter space is a product space. So, for example, we know which $h$ would maximize the likelihood of the data, without having to know anything about the coarsening mechanism (except that it's CAR, of course). It of course also implies lots more good consequences for likelihood-based (and, in particular, Bayesian) inference.

PROOF OF PROPOSITION 1.2. Let $k$ be a positive measurable function on $\mathcal{X}$. Remember that

$$\mu(dz) = g \circ \psi(z) \cdot h \circ \pi(z)\mu_0(dz).$$

Then we have

$$\begin{aligned} E_\mu(k(X)) &= E_{\mu_0}(k(X)g(X)h(Y)) \\ &= E_{P_0}(k(X)g(X)E_{\mu_0}(h(Y)|X)) \\ &= E_{P_0}(k(X)g(X)S(h)(X)). \quad \square \end{aligned}$$

The CAR assumption as we defined it depends on the choice of $\mu_0$, but we do have the following proposition:

PROPOSITION 1.3. *Let $\mu_0$ and $\nu_0$ be probability measures on $\mathcal{Z}$ such that $\nu_0$ satisfies the CAR assumption for $\mu_0$ (in particular, $\nu_0 \ll \mu_0$). Then a probability measure $\mu \ll \nu_0$ on $\mathcal{Z}$ satisfies the CAR assumption for $\mu_0$ if and only if it satisfies the CAR assumption for $\nu_0$.*

PROOF. Since $\nu_0$ satisfies the CAR assumption for $\mu_0$, we can write
$$\nu_0(dz) = g_0 \circ \psi(z) h_0 \circ \pi(z) \mu_0(dz)$$
such that $h_0$ is a density for $Q_0$ and $E_{\mu_0}(g_0(X)|Y) = 1$, which means that $\nu_0(dz|y) = g_0 \circ \psi(z) \mu_0(dz|y)$. Suppose $\mu$ satisfies CAR for $\mu_0$, so we can write
$$\mu(dz) = g_1 \circ \psi(z) h_1 \circ \pi(z) \mu_0(dz)$$
with $E_{\mu_0}(h_1(Y)) = 1$ and $E_{\mu_0}(g_1(X)|Y) = 1$. Note that
$$h_1 Q_0 = \pi(\mu) \ll \pi(\nu_0) = h_0 Q_0,$$
so $h_1/h_0$ is well defined ($0/0 = 0$). The same reasoning, but with $\pi$ replaced with $\psi$, gives that $g_1/g_0$ is well defined. Now note that
$$E_{\nu_0}\left(\frac{h_1}{h_0}(Y)\right) = E_{\mu_0}(h_1(Y)) = 1$$
and
$$E_{\nu_0}\left(\frac{g_1}{g_0}(X)\Big|Y=y\right) = \int_{\mathcal{Z}} \frac{g_1}{g_0}(\psi(z)) \nu_0(dz|y)$$
$$= \int_{\mathcal{Z}} \frac{g_1}{g_0}(\psi(z)) g_0(\psi(z)) \mu_0(dz|y)$$
$$= \int_{\mathcal{Z}} g_1(\psi(z)) \mu_0(dz|y)$$
$$= 1,$$
so
$$\mu(dz) = (g_1/g_0) \circ \psi(z)(h_1/h_0) \circ \pi(z) \nu_0(dz)$$
satisfies CAR for $\nu_0$.

If $\mu$ satisfies CAR for $\nu_0$, we conclude in a completely analogous way that $\mu$ satisfies CAR for $\mu_0$. □

This proposition shows that for any $\mu_0$ you pick such that a certain coarsening mechanism $\nu_0$ satisfies CAR for $\mu_0$ (and is, therefore, an element of your model), the possible distributions of $Z$ absolutely continuous with respect to $\nu_0$ are the same as when you would have chosen $\mu_0 = \nu_0$. Therefore, a logical choice for $\mu_0$ is a generic distribution for $Z$ that you would want to have in your model, preferably with an as large as possible support.

One can easily verify that our definition of the CAR assumption is equivalent to the ones given in Gill, van der Laan and Robins (1997) (for the dominated case), Jacobsen and Keiding (1995) and Nielsen (2000), when we restrict ourselves to their respective setups (see also the discussion after Theorem 3.8). We would like to point out that for the factorization property of Proposition 1.2, Gill, van der Laan and Robins (1997) also have to restrict themselves to the dominated case. The conjecture made in Gill, van der Laan and Robins (1997) is that the CAR assumption does not restrict the possible distributions of the data, making it impossible to test whether the CAR assumption is fulfilled or not. In fact, they prove this conjecture (in their setup) when $\mathcal{Y}$ is a finite space. In the next section we will give examples where the conjecture actually *fails*, not only in our generalized setup, but also in the more restrictive setups. In Section 3 we will give sufficient and almost necessary conditions when the conjecture will hold.

## 2. Examples.

EXAMPLE 2.1. Let $\mathcal{Y} = [0, \infty[$, $\mathcal{Z} = [0, \infty[ \times [0, \infty[$ and $Z = (Y, C)$. Define $X = \psi(Y, C) = CY$. This coarsening mechanism cannot be described as in Gill, van der Laan and Robins (1997), for knowing $X$ is not equivalent to knowing that $Y$ lies in the set of points compatible with the observation $X$. Now we have to choose $\mu_0$:

$$\mu_0(dy\,dc) = e^{-y}e^{-c}\,dy\,dc.$$

The CAR assumption states that for a possible distribution $\mu$ of $Z$, there exist $h \in L^1(Q_0)$ and $g \in L^1(P_0)$ such that

$$\mu(dy\,dc) = (g(cy)e^{-c}\,dc)h(y)e^{-y}\,dy.$$

Furthermore, (1.1) tells us that

$$\int_0^\infty g(cy)e^{-c}\,dc = 1 \qquad \forall y > 0.$$

But this means that the Laplace transform of $g$ is identically equal to the Laplace transform of 1, and, therefore, $g = 1$. So the possible choices for $\mu$ are

$$\mu(dy\,dc) = h(y)e^{-y}e^{-c}\,dy\,dc,$$

where $h$ is a density with respect to $Q_0(dy) = e^{-y}\,dy$. Note that $C$ is independent of $Y$ with a given distribution, and the distribution of $Y$ is arbitrary. A simple transformation of variables gives

$$\psi(\mu)(dx) = \left(\int_0^\infty h(y)e^{-y}e^{-x/y}\frac{1}{y}\,dy\right)dx.$$

Therefore, $X$ always has a decreasing density with respect to the Lebesgue measure on $[0, \infty[$, which shows that in this case the CAR assumption does restrict the possible distributions of the data.

As noted before, the CAR assumption depends on the choice of $\mu_0$. To illustrate this, let us choose

$$\mu_0(dy\,dc) = (ye^{-yc}\,dc) \cdot e^{-y}\,dy.$$

Then CAR implies for our (positive) function $g$ that

$$\int_0^\infty g(cy)ye^{-cy}\,dc = 1 \qquad \forall y > 0.$$

However, this is nothing more than saying that $g$ is a density for $P_0$, since in this case $P_0$ is the standard exponential! Clearly, this means that the CAR assumption is not testable in this case. However, it is not hard to see that in this case $S(h) = 1$, so all information about $Y$ is lost. As a final remark, note that the CAR assumption is only affected by $\mu_0$ through $\mu_0(dc|y)$, the conditional distribution of $C$ given $Y = y$, so that choosing a different (but equivalent) $Q_0$ essentially leaves the CAR assumption unaltered (this also follows from Proposition 1.3).

EXAMPLE 2.2 (Current status). A much more important example, and one that also fits the setups of Jacobsen and Keiding (1995) and Gill, van der Laan and Robins (1997), is that of current status data. We will consider the bounded case, that is, all times considered fall in $[0, 1]$, but it is not hard to see that this is not a real restriction. So define $\mathcal{Y} = [0, 1]$, $Y$ is the time of interest, $C \in [0, 1]$, the censoring time, and $Z = (Y, C)$, so $\mathcal{Z} = [0, 1] \times [0, 1]$. Define

$$\psi(Y, C) = (C, \mathbb{1}_{\{Y \leq C\}}),$$

so $\mathcal{X} = [0, 1] \times \{0, 1\}$. The interpretation is that one knows the time one visited the doctor, and the doctor can say whether someone is sick or not. Choose $\mu_0(dy\,dc) = dy\,dc$. Then (1.1) implies that we can choose positive $g \in L^1(P_0)$ such that

$$\int_0^1 g(c, \mathbb{1}_{\{y \leq c\}})\,dc = 1 \qquad \forall 0 \leq y \leq 1.$$

However, this says that

$$\int_0^y g(c, 0)\,dc + \int_y^1 g(c, 1)\,dc = 1 \qquad \forall 0 \leq y \leq 1.$$

Differentiating with respect to $y$ shows that

$$g(c, 0) = g(c, 1) \qquad \forall 0 \leq c \leq 1.$$

So CAR implies that the only allowed models for $\mu$ are

$$\mu(dy\,dc) = g(c)h(y)\,dc\,dy,$$

where $h$ and $g$ are densities on $[0,1]$. This is, of course, equivalent with saying that $Y$ and $C$ have to be independent.

Consider the following subsets of $\mathcal{X}$:

$$A_1 = \{(x,1) : x \in [0,\tfrac{1}{2}]\} \quad \text{and} \quad A_2 = \{(x,0) : x \in [\tfrac{1}{2},1]\}.$$

Let $\mathcal{P}$ be the set of all probability distributions on $\mathcal{X}$ and define for every $P \in \mathcal{P}$,

$$\Phi(P) = (P(A_1), P(A_2)).$$

Clearly,

$$\Phi(\mathcal{P}) = \{(a_1,a_2) \in [0,1]^2 : a_1 + a_2 \leq 1\}.$$

Now suppose the CAR assumption holds, so $Y$ and $C$ are independent. Then we know that

$$\mathbb{P}(X \in A_1) = \mathbb{P}(C \leq \tfrac{1}{2} \text{ and } Y \leq C)$$
$$\leq \mathbb{P}(C \leq \tfrac{1}{2} \text{ and } Y \leq \tfrac{1}{2})$$
$$= \mathbb{P}(C \leq \tfrac{1}{2}) \cdot \mathbb{P}(Y \leq \tfrac{1}{2}).$$

Similarly,

$$\mathbb{P}(X \in A_2) \leq \mathbb{P}(C \geq \tfrac{1}{2}) \cdot \mathbb{P}(Y \geq \tfrac{1}{2}).$$

This means that

$$P(X \in A_1) \cdot P(X \in A_2) \leq \tfrac{1}{16}.$$

So, if we define $\mathcal{P}_{\text{CAR}}$ as the set of all possible distributions of the data under the CAR assumption, then

$$\Phi(\mathcal{P}_{\text{CAR}}) \subset \{(a_1,a_2) \in [0,1]^2 : a_1 + a_2 \leq 1 \text{ and } a_1 \cdot a_2 \leq \tfrac{1}{16}\}.$$

Since this is a proper subset of $\Phi(\mathcal{P})$, we conclude that in the case of current status, it is possible to find a distribution of the data that contradicts the CAR assumption. In a future paper we will discuss what would be a good way to test the CAR assumption in this important example. Here we would like to note a few things. In the first place, it is possible that the data distribution is an element of $\mathcal{P}_{\text{CAR}}$, even though the CAR assumption is not fulfilled: one easily checks that this happens when

$$c \mapsto \int_0^c f(y|c)\,dy$$

is a continuous distribution function (i.e., nondecreasing), where $f(y|c)$ is the conditional density of $Y$ given $C = c$. This shows that it is impossible to verify CAR by the data; it is just sometimes possible to reject the CAR assumption.

In the second place we note that the convex hull of all independent densities of $(Y,C)$ is weakly dense in the set of all densities, and, therefore, the convex hull of $\mathcal{P}_{\text{CAR}}$ is weakly dense in $\mathcal{P}$. This means that you cannot test the CAR assumption with one linear test function. In particular, it shows that the model for the distribution of the data under CAR is not convex.

As a third remark, we would like to point out to the reader that although this example fits in the setup of Gill, van der Laan and Robins (1997) for CAR on general sample spaces, it does not fit in their setup for finite spaces, not even when we restrict $Y$ and $C$ to finitely many possible values. This is because the observed sets are all of the form $\{Y \leq C\}$ or $\{Y > C\}$, and it is essential in their setup that the CAR assumption allow distributions such that all possible nonempty subsets of $\mathcal{Y}$ might be observable. See also the discussion after Theorem 3.8.

Finally, it is not hard to show that under the assumption CAR(ABS) defined in Gill, van der Laan and Robins (1997), one can find all possible distributions of $X$ by assuming that $Y$ and $C$ are independent, but can have any distribution (not necessarily dominated). This means that the argument given here also shows that CAR(ABS) restricts the possible distributions of the data $X$. We do not think that by restricting ourselves to the dominated case we throw away an important part of the possible distributions of $X$ under CAR(ABS).

**3. General conditions for the testability of CAR.** In this section we will give our most abstract definition of coarsened data, but we will first look at the map $S: L^1(Q_0) \to L^1(P_0)$. We will repeat its definition:

(3.1) $$S(h)(x) = E_{\mu_0}(h(Y)|X = x).$$

If we denote the duality between $L^1$-functions and $L^\infty$-functions by $\langle \cdot, \cdot \rangle$, we would like to remind the reader that the dual map

$$S^*: L^\infty(P_0) \to L^\infty(Q_0)$$

is defined such that

$$\langle S(h), k \rangle = \langle h, S^*(k) \rangle.$$

Note that for $k \in L^\infty(P_0)$,

$$\langle S(h), k \rangle = E_{\mu_0}(k(X)E_{\mu_0}(h(Y)|X)) = E_{\mu_0}(k(X)h(Y)).$$

PROPOSITION 3.1. *The linear map $S: L^1(Q_0) \to L^1(P_0)$ defined above has the following properties:*

1. $S(1) = 1$ and $S^*(1) = 1$, where $S^*$ denotes the dual of $S$.
2. $S$ is positive, that is, $h \geq 0 \Rightarrow S(h) \geq 0$.
3. $\|S\| = 1$, where $\|\cdot\|$ denotes the operator-norm.

PROOF.    Properties 1 and 2 are obvious. It is also clear that
$$\|h \circ \pi\|_1 = \|h\|_1 \quad \text{and} \quad \|k \circ \psi\|_\infty = \|k\|_\infty$$
[here we use $Q_0 = \pi(\mu_0)$ and $P_0 = \psi(\mu_0)$], which shows that $\|S\| \leq 1$. Since $S(1) = 1$, $\|S\| = 1$.  □

The importance of the map $S$ is seen most clearly when we translate (1.1):
$$E_{\mu_0}(g(X)|Y) = 1.$$
It is well known that
$$S^*(g)(y) = E_{\mu_0}(g(X)|Y = y),$$
so this means that the CAR assumption restricts our choice for $g$ (remember that $g \circ \psi$ is the conditional density of $Z$ given $Y = y$, for all $y$) to all positive $g$ such that
$$S^*(g) = 1.$$
This will lead us to a new definition of CAR.

DEFINITION 3.2.    Let $Y$ be a stochastic variable of interest, defined on a space $\mathcal{Y}$, and let $Q_0$ be a probability measure on $\mathcal{Y}$. Let $\mathcal{X}$ be the data-space and $P_0$ a probability measure on $\mathcal{X}$. We define a *coarsening* (of $Y$) as a linear map
$$S : L^1(Q_0) \to L^1(P_0)$$
such that:

1. $S(1) = 1$ and $S^*(1) = 1$, where $S^*$ denotes the dual of $S$.
2. $S$ is positive, that is, $h \geq 0 \Rightarrow S(h) \geq 0$.

We thank one of the referees for pointing out the following: every coarsening $S$ can be obtained through a conditional expectation, as we did in the original definition of CAR. To see this, define $\mathcal{Z} = \mathcal{Y} \times \mathcal{X}$. We define a probability measure $\mu_0$ on $\mathcal{Y} \times \mathcal{X}$ in the following way: let $A \subset \mathcal{Y}$ and $B \subset \mathcal{X}$ be measurable such that $\mathbb{1}_A \in L^1(Q_0)$ and $\mathbb{1}_B \in L^1(P_0)$. Then we define
$$\mu_0(A \times B) = E_{P_0}(\mathbb{1}_B(X) S(\mathbb{1}_A)(X)).$$
This extends to a probability measure on $\mathcal{Y} \times \mathcal{X}$ such that for $h \in L^1(Q_0)$ and $k \in L^1(P_0)$,
$$E_{\mu_0}(k(X)h(Y)) = E_{P_0}(k(X)S(h)(X)).$$

It is easy to check that $Q_0$ and $P_0$ are the marginals of $Y$, respectively, $X$, and that

$$S(h)(x) = E_{\mu_0}(h(Y)|X = x).$$

From this it is clear that

$$S^*(k)(y) = E_{\mu_0}(k(X)|Y = y),$$

so the map $S^*$ is in itself a coarsening of $X$. This is the content of the next lemma, which we will prove without using the auxiliary measure $\mu_0$. In fact, we believe the map $S$ to be the most convenient object to study, which is why we will not refer to $\mu_0$ again.

LEMMA 3.3. *Let $S: L^1(Q_0) \to L^1(P_0)$ be a coarsening. Then:*

1. $S$ *is continuous and* $\|S\| = 1$.
2. *The dual map $S^*$ is also defined and continuous from $L^1(P_0)$ to $L^1(Q_0)$ (in fact, $S^*$ is a coarsening itself).*

PROOF. Let $h \in L^1(Q_0)$. Then $-|h| \leq h \leq |h|$, so $|S(h)| \leq S(|h|)$. Now,

$$\|S(|h|)\| = \langle S(|h|), 1 \rangle$$
$$= \langle |h|, 1 \rangle$$
$$= \|h\|.$$

This, together with $S(1) = 1$, proves the first statement.

Let $g \in L^1(P_0)_+$. There exists $\{g_n\} \subset L^\infty(P_0)_+$ such that $g_n \uparrow g$. Clearly, $S^*$ is also positive, so $S^*(g_n) \uparrow h$ for some $h \in L^1(Q_0)$ [note that $\langle h, 1 \rangle = \lim \uparrow \langle S^*(g_n), 1 \rangle = \lim \uparrow \langle g_n, S(1) \rangle = \langle g, 1 \rangle$]. Also, if $h' \in L^\infty(Q_0)$, then $\langle h, h' \rangle = \langle g, S(h') \rangle$, so $h$ does not depend on the sequence $\{g_n\}$. Define $h = S^*(g)$. It is trivial to check that with this definition, $S^*$ is in itself a coarsening. $\square$

Define for a probability measure $\nu$,

$$\Delta_\nu = \{h \in L^1(\nu) : h \geq 0 \text{ and } \langle h, 1 \rangle = 1\},$$

the set of densities with respect to $\nu$.

DEFINITION 3.4 (CAR). Let

$$S: L^1(Q_0) \to L^1(P_0)$$

be a coarsening of a random variable $Y$. The CAR assumption now states that the distribution of the data belongs to the set

$$\mathcal{P}_{\text{CAR}} = \{g \cdot S(h) : h \in \Delta_{Q_0}, g \in L^1(P_0)_+ \text{ and } S^*(g) = 1\}.$$

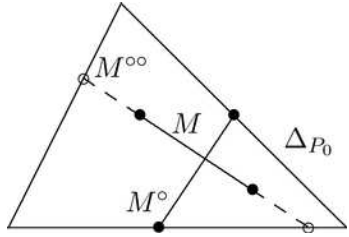

Fig. 1.

First we should note that $\mathcal{P}_{\mathrm{CAR}} \subset \Delta_{P_0}$, because
$$\langle g \cdot S(h), 1 \rangle = \langle h, S^*(g) \rangle = 1$$
and $S$ is a positive map. In this new definition we also retain the product structure of the likelihood of the data. The remark after Definition 3.2 shows that the only difference with the previous definition is that we restrict the distributions of the data $X$, instead of restricting the distributions of the hidden variable $Z$.

It is clear that the question of testability of the CAR assumption amounts to checking whether the set $\mathcal{P}_{\mathrm{CAR}}$ is dense in $\Delta_{P_0}$. Before we consider this question, we want to note the following: define
$$M = \{S(h) \in \Delta_{P_0} : h \in \Delta_{Q_0}\}.$$
Then $M$ is a convex subset of $\Delta_{P_0}$. Now in analogy to the polar set of a subspace of a linear space, we define
$$M^\circ = \{g \in L^1(P_0)_+ : (\forall\, h \in M) \langle h, g \rangle = 1\}.$$
Since for all $g \in L^1(P_0)_+$, $S^*(g) = 1$ is equivalent to
$$(\forall\, h \in \Delta_{Q_0}) \qquad \langle S(h), g \rangle = 1,$$
we get that
$$M \cdot M^\circ = \mathcal{P}_{\mathrm{CAR}}.$$
Encouraged by this observation, we define
$$M^{\circ\circ} = \{h \in L^1(P_0)_+ : (\forall\, g \in M^\circ) \langle h, g \rangle = 1\}.$$

Figure 1 shows the situation when $P_0$ has a support of 3 points (so we can view $\Delta_{P_0}$ as a triangle) and $M$ is a convex subset of $\Delta_{P_0}$.

As you can see, we should view $M^{\circ\circ}$ as an extension of $M$ to the edges of $\Delta_{P_0}$. The following proposition, together with Lemma 3.9, substantiates Figure 1:

PROPOSITION 3.5. *Let $M$ be an arbitrary subset of $\Delta_P$, with $P$ some probability measure.*

1. $M \subset M^{\circ\circ} \subset \Delta_P$.
2. $(M^{\circ\circ})^\circ = M^\circ$.

PROOF.  1. The first inclusion is obvious from the definitions. For the second one, it is enough to note that $1 \in M^\circ$, because $M \subset \Delta_P$.

2. Clearly, $(M^{\circ\circ})^\circ \subset M^\circ$. Let $g \in M^\circ$. If $h \in (M^\circ)^\circ$, then $\langle h, g \rangle = 1$ (because $g \in M^\circ$), so $g \in (M^{\circ\circ})^\circ$.  □

Since $M \cdot M^\circ \subset M^{\circ\circ} \cdot M^\circ \subset \Delta_P$, a natural necessary condition on $M$ for $M \cdot M^\circ$ to be dense in $\Delta_P$ would be $M^{\circ\circ} \subset \overline{M}$. The following proposition more or less substantiates this statement. We do have to caution the reader that in principle $M^\circ$ and $M^{\circ\circ}$ need not be closed sets, since the linear functional $h \mapsto \langle h, g \rangle$ is not continuous on $L^1(P)$ if $g \in L^1(P)_+ \setminus L^\infty(P)$.

PROPOSITION 3.6.  *Let $M$ be a subset of $\Delta_P$ such that $M^{\circ\circ} \not\subset \overline{M}$. Then there exist $h \in M^{\circ\circ}$ and $\varepsilon > 0$ such that for all $f \in M \cdot M^\circ$,*

$$\int -\log\left(\frac{f}{h}\right) h \, dP > \varepsilon.$$

PROOF.  Choose $h \in M^{\circ\circ} \setminus \overline{M}$. Then there exists $\varepsilon > 0$ such that for all $\tilde{h} \in M$, $\|h - \tilde{h}\| > \varepsilon$. It is a well-known inequality for the Kullback–Leibler divergence [see, e.g., van der Vaart (1998), page 62] that

$$\int -\log\left(\frac{\tilde{h}}{h}\right) h \, dP \geq \frac{1}{4} \|h - \tilde{h}\|^2.$$

Now let $f \in M \cdot M^\circ$, so $f = \tilde{h} g$, with $\tilde{h} \in M$ and $g \in M^\circ$. Note that $\langle g, h \rangle = 1$, since $h \in M^{\circ\circ}$. So

$$\int -\log\left(\frac{\tilde{h}g}{h}\right) h \, dP = \int -\log\left(\frac{\tilde{h}}{h}\right) h \, dP + \int -\log(g) h \, dP$$
$$> \frac{1}{4}\varepsilon^2 - \log\left(\int gh \, dP\right)$$
$$= \frac{1}{4}\varepsilon^2. \quad □$$

We have to point out that this proposition does not state that, under the assumption that $M^{\circ\circ} \not\subset \overline{M}$, $M \cdot M^\circ$ is *not* dense in $\Delta_P$. We were not able to prove that statement in general. However, it does *indicate* that $M \cdot M^\circ$ is not dense in $\Delta_P$, and in specific examples it should not be too hard to actually prove it.

EXAMPLE 3.7 (Current status). As we have seen already, we consider a time of interest $Y \in [0,1]$, a censoring time $C \in [0,1]$, and the data consists of $(C, \mathbb{1}_{\{Y \leq C\}})$. We take

$$Q_0(dt) = dt \quad \text{and} \quad P_0(dx, \delta) = x\,dx \cdot \mathbb{1}_{\{\delta=1\}} + (1-x)\,dx \cdot \mathbb{1}_{\{\delta=0\}}.$$

It is easily seen that our map $S$ is equal to

$$S(h)(x,\delta) = \frac{1}{x}\int_0^x h(t)\,dt \cdot \mathbb{1}_{\{\delta=1\}} + \frac{1}{1-x}\int_x^1 h(t)\,dt \cdot \mathbb{1}_{\{\delta=0\}}.$$

Remember that $M = S(\Delta_{Q_0})$, so for all $\tilde{h} \in \overline{M}$, we have that $x\tilde{h}(x,1)$ is increasing in $x$. Now choose

$$h(t) = \mathbb{1}_{\{t \leq 1/3\}} - \mathbb{1}_{\{1/3 < t \leq 2/3\}} + 3 \cdot \mathbb{1}_{\{t > 2/3\}}.$$

Then $\langle h, 1 \rangle = 1$ and $S(h) \geq 0$, so $S(h) \in M^{\circ\circ}$, but $xS(h)(x,1) = \int_0^x h(t)\,dt$ is not increasing in $x$, so $S(h) \notin \overline{M}$. It was this observation that led us to find the test described in Example 2.2.

The statement we would like to prove for $M \subset \Delta_P$ is that $M \cdot M^\circ$ is dense in $\Delta_P$ if and only if $M^{\circ\circ} \subset \overline{M}$. However, we were not able to prove it in this generality, nor find a counterexample to it. Only when $P$ has finite support were we able to prove the statement in full generality:

THEOREM 3.8. *Let $P$ be a probability measure with finite support and let $M \subset \Delta_P$ such that there exists $h_0 \in M$ with $h_0 > 0$. Then $M \cdot M^\circ$ is dense in $\Delta_P$ if and only if $\overline{M} = M^{\circ\circ}$.*

PROOF. Let $\overline{M} \neq M^{\circ\circ}$. Since we are now in the situation that $L^1(P) = L^\infty(P)$, it follows that $M^{\circ\circ}$ is closed, so we always have $\overline{M} \subset M^{\circ\circ}$. According to Proposition 3.6, there exist $h \in M^{\circ\circ}$ and $\varepsilon > 0$ such that for all $f \in M \cdot M^\circ$,

$$\int -\log\left(\frac{f}{h}\right) h\,dP > \varepsilon.$$

Since $h_0 \in M^{\circ\circ}$, we can choose $h > 0$ [note that $\varepsilon h_0 + (1-\varepsilon)h \in M^{\circ\circ}$, for all $1 > \varepsilon > 0$]. Since $\{f > 0 : f \in \Delta_P\}$ is an open subset of $\Delta_P$ and since $f \mapsto \int -\log(f/h)h\,dP$ is continuous on this set (so, in particular, continuous at $h$), we conclude that there exists $\eta > 0$ such that for all $f \in M \cdot M^\circ$, $\|f - h\| > \eta$.

Now let $\overline{M} = M^{\circ\circ}$. Choose $f \in \Delta_P$ with $f > 0$. Since $\overline{M}$ is compact and $h \mapsto \int -\log(h/f)f\,d\mu$ is lower semi-continuous (see also Lemma 3.11), there exists $h \in \overline{M}$ that minimizes this Kullback–Leibler divergence. It is also clear that $h > 0$, since otherwise the Kullback–Leibler divergence would be $+\infty$ (here we use that $h_0 \in M$). Now let $\tilde{h} \in M$. Since $h > 0$, there exists $\varepsilon > 0$ such that when $|\lambda| < \varepsilon$, $h + \lambda(\tilde{h} - h) \geq 0$. This means that $h + \lambda(\tilde{h} - h) \in M^{\circ\circ} = \overline{M}$, because clearly $\langle h + \lambda(\tilde{h} - h), g \rangle = 1$ for all $g \in M^\circ$. The function

$$\lambda \mapsto \int -\log\left(\frac{h + \lambda(\tilde{h} - h)}{f}\right) f \, dP$$

has a minimum at $\lambda = 0$ for $\lambda \in \,]-\varepsilon, -\varepsilon[$, so the derivative at $\lambda = 0$ (which exists!) must be zero. A simple calculation yields

$$\int (\tilde{h} - h) \frac{f}{h} \, dP = 0.$$

This proves that $\langle \tilde{h}, f/h \rangle = 1$ for all $\tilde{h} \in M$, so $f/h \in M^\circ$. Therefore, $f \in \overline{M} \cdot M^\circ$. It is not hard to see that if $h_n \to h$, then $h_n \cdot f/h \to f$, which proves that $M \cdot M^\circ$ is dense in $\Delta_P$. $\square$

This theorem is very much like the theorem in Section 2 of Gill, van der Laan and Robins (1997) and also the proof is very similar. To show how their theorem (apart from the uniqueness statement) follows from Theorem 3.8, we translate their setup into ours. Let $\mathcal{Y}$ be a finite space with $m$ points and let $\mathcal{X} = \mathcal{P}(\mathcal{Y}) \setminus \{\varnothing\}$, the collection of all nonempty subsets of $\mathcal{Y}$. The idea is that one observes $X \subset \mathcal{Y}$ such that $Y \in X$. To reformulate the CAR assumption used in Gill, van der Laan and Robins (1997), we define $\mathcal{Z} = \{(y, A) : y \in A \subset \mathcal{Y}\}$ and $\mu_0$ as the rescaled counting measure on $\mathcal{Z}$, assigning mass $2^{1-m}/m$ to each element of $\mathcal{Z}$. Obviously, we define $\pi(y, A) = y$ and $\psi(y, A) = A$, so $Q_0 = \pi(\mu_0)$ is the rescaled counting measure on $\mathcal{Y}$ (assigning mass $1/m$ to each point) and $P_0 = \psi(\mu_0)$ satisfies

$$P_0(\{A\}) = \frac{|A|}{m 2^{m-1}} \qquad (\forall A \subset \mathcal{Y}),$$

where $|A|$ denotes the number of elements of $A$. Now we define $S \colon L^1(Q_0) \to L^1(P_0)$ such that for all $h \in L^1(Q_0)$ and $A \in \mathcal{X}$, we have

$$(3.2) \qquad S(h)(A) = E_{\mu_0}(h(Y)|X = A) = \frac{1}{|A|} \sum_{y \in A} h(y).$$

It follows immediately that for $g \in L^1(P_0)$ and $y \in \mathcal{Y}$, we have

$$S^*(g)(y) = 2^{1-m} \sum_{A \ni y} g(A).$$

The CAR assumption now states that the likelihood of $X$ with respect to $P_0$ equals $g \cdot S(h)$, where $h$ is an arbitrary density with respect to $Q_0$ and $g \in L^1(P_0)_+$ such that $S^*(g) = 1$. If we would follow Definition 1.1, we would restrict the possible distributions $\mu$ of $Z = (Y, X)$ such that

$$\mu(X = A | Y = y) = g(A) \mu_0(X = A | Y = y) = 2^{1-m} g(A) \mathbb{1}_{\{y \in A\}}.$$

It is not hard to see that this is indeed equivalent to the definition of Gill, van der Laan and Robins (1997) used for finite sample spaces. So, in fact, they use a very specific form of the map $S$; even in finite sample spaces our setup is much less restrictive. Finally, to conclude that in this case CAR cannot be tested, we use Theorem 3.8 to see that we only need to check that when we define

$$M = \{S(h) : h \in \Delta_{Q_0}\},$$

we have $M^{\circ\circ} = \overline{M}$. We will use the following lemma.

LEMMA 3.9. *Let $P$ be a measure with finite support and let $M \subset \Delta_P$. Then*

$$M^{\circ\circ} = \langle M \rangle \cap \Delta_P.$$

*Here $\langle M \rangle$ denotes the linear span of $M$.*

PROOF. Let $h \in \langle M \rangle \cap \Delta_P$, so $h = \sum \lambda_i h_i$ with $\lambda_i \in \mathbb{R}$ and $h_i \in M$ such that $h \geq 0$ and $\langle h, 1 \rangle = 1$. This means that $\sum \lambda_i = 1$. If $g \in M^{\circ}$, then for every $i$ $\langle h_i, g \rangle = 1$, so we conclude that $\langle h, g \rangle = 1$, and, therefore, $h \in M^{\circ\circ}$. We have shown that $\langle M \rangle \cap \Delta_P \subset M^{\circ\circ}$.

Now suppose $h \in \Delta_P$ and $h \notin \langle M \rangle$. Since $L^1(P)$ is finite dimensional, there exists $\phi \in L^1(P)$ such that for all $\tilde{h} \in \langle M \rangle$, we have $\langle \tilde{h}, \phi \rangle = 0$ and $\langle h, \phi \rangle > 0$. We can choose $\phi$ such that $|\phi| \leq 1$. Define $g = 1 + \phi$. Then $g \geq 0$ and for $\tilde{h} \in M$ we have $\langle \tilde{h}, g \rangle = 1$, so $g \in M^{\circ}$. However, $\langle h, g \rangle > 1$, so $h \notin M^{\circ\circ}$. Since $M^{\circ\circ} \subset \Delta_P$, we have shown that $M^{\circ\circ} \subset \langle M \rangle \cap \Delta_P$. □

When $M = \{S(h) : h \in \Delta_{Q_0}\}$, it is easy to check that $\langle M \rangle \cap \Delta_{P_0} = \{S(h) : h \in L^1(Q_0), \langle h, 1 \rangle = 1, S(h) \geq 0\}$. Therefore, whenever $\mathcal{X}$ is a finite set, $M^{\circ\circ} = \overline{M}$ is equivalent to

(3.3) $\quad S(h) \geq 0 \quad \implies \quad \exists \tilde{h} \geq 0 : S(\tilde{h}) = S(h) \qquad [\forall h \in L^1(Q_0)].$

For the map $S$ we were considering, this follows trivially from (3.2) [note that $S(h)(\{y\}) = h(y)$].

The problem with extending the proof of Theorem 3.8 to general $P$ is twofold. First of all, $\overline{M}$ will not be compact in general, which makes it difficult to find a minimum for the Kullback–Leibler divergence. The second problem is concluding that the derivative is zero: even if we find a minimum (in some compactification), we can only conclude that the directional derivative we used in the previous proof is negative, but not necessarily zero. To solve these problems and come up with a theorem that can be used for practical situations, we will use the map $S$ more extensively by putting restrictions on it. But first we will discuss an extension of the Kullback–Leibler divergence to solve the noncompactness problem.

DEFINITION 3.10. Let $E = (L^\infty(P))'$, the (strong) dual of the Banach space $L^\infty(P)$. Let $f \in L^1(P)_+$. Define for $h \in E$, $h \geq 0$,

$$\mathrm{KL}_f(h) = \sup\left\{\sum_{i=1}^n -\log\left(\frac{\langle h, \phi_i\rangle}{\langle f, \phi_i\rangle}\right)\langle f, \phi_i\rangle : \phi_i \in L^\infty(P)_+, \sum_{i=1}^n \phi_i = 1\right\}.$$

We would like to make a few remarks. As $E$ is the dual of an ordered Banach space, it is itself ordered in the obvious way: $h \geq 0$ if for all $\phi \in L^\infty(P)_+$, $\langle h, \phi\rangle \geq 0$. Furthermore, $L^1(P) \subset E$. We also have that the unit ball of $E$ is weakly compact (Banach–Alaoglu), and if $h \in E_+$ (i.e., $h$ is positive), we have that $\|h\| = \langle h, 1\rangle$. Since $\mathrm{KL}_f$ is the supremum of weakly continuous functions on $E_+$, it is itself weakly lower semi-continuous. If $M \subset \Delta_P$, then $\overline{M}^\sigma$ (the closure of $M$ in the weak topology, seen as a subset of $E$) will be weakly compact, because $\overline{M}^\sigma \subset E_+$ and for all $h \in \overline{M}^\sigma$, $\langle h, 1\rangle = 1$, so it is a weakly closed subset of the unit ball. This means that $\mathrm{KL}_f$ will attain its minimum on $\overline{M}^\sigma$ for some $h \in \overline{M}^\sigma$.

From the theory of ordered vector lattices [see, e.g., Schaefer and Wolff (1999), Chapter V] it follows that $L^1(P)$ is a *band* in $E$. This means that each $h \in E_+$ can be uniquely decomposed as $h = h_{/\!/} + h_\perp$, where $h_{/\!/} \in L^1(P)_+$ and $h_\perp \geq 0$ is disjoint from $L^1(P)$, so for all $f \in L^1(P)_+$, we have that $\inf(h_\perp, f) = 0$ (compare this to the decomposition of a measure into a part that is absolutely continuous to some other measure and a part which is disjoint from this other measure). We have the following lemma, the proof of which is deferred to the Appendix.

LEMMA 3.11. *Let $f \in L^1(P)_+$. Then, in the notation introduced above, for all $h \in E_+$,*

$$\mathrm{KL}_f(h) = \mathrm{KL}_f(h_{/\!/}) = \int -\log\left(\frac{h_{/\!/}}{f}\right) f\, dP.$$

Now we will consider a coarsening $S: L^1(Q_0) \to L^1(P_0)$. Define $E_{Q_0} = (L^\infty(Q_0))'$ and $E_{P_0} = (L^\infty(P_0))'$. By considering the dual map of $S^*$, we can extend $S: E_{Q_0} \to E_{P_0}$. Clearly, $S$ will be continuous for the weak topologies on $E_{Q_0}$ and $E_{P_0}$ (as well as for the strong topologies) and $S$ will be a positive map. Define $M = S(\Delta_{Q_0})$. Since $\overline{\Delta_{Q_0}}^\sigma \subset E_{Q_0}$ is weakly compact, $\overline{M}^\sigma = S(\overline{\Delta_{Q_0}}^\sigma)(\subset E_{P_0})$. When $h \in E_{Q_0,+}$, we can consider $h_{/\!/} \in L^1(Q_0)_+$ as well as $S(h)_{/\!/} \in L^1(P_0)_+$. In general, we can only deduce that $S(h)_{/\!/} \geq S(h_{/\!/})$, since $h = h_{/\!/} + h_\perp$ and $S(h_{/\!/}) \in L^1(P_0)_+$.

Before we can state our main result, we need two assumptions. The first is the analogue of $\overline{M} = M^{\circ\circ}$, or equation (3.3) which we discussed before, but slightly stronger:

(A1) For all $h' \in E_{Q_0,+}$ such that $S(h')_{/\!/} > 0$, there exists $h \in E_{Q_0,+}$ with $S(h) = S(h')$ and $h_{/\!/} > 0$.

How we will use assumption (A1) is stated in the following lemma: we say that $h_1 \in L^1(P)_+$ *dominates* $h_2 \in L^1(P)_+$ (notation: $h_2 \lesssim h_1$), if there exists $R > 0$ such that $h_2 \leq R h_1$.

LEMMA 3.12. *Suppose $h_0 \in E_{Q_0,+}$ such that $h_{0,/\!/} > 0$. Let $h \in L^1(Q_0)_+$. Then there exists a sequence $h_n \in L^1(Q_0)_+$ such that $h_n \lesssim h_{0,/\!/}$ and $h_n \uparrow h$.*

PROOF. Define $f = h_{0,/\!/}$. Let $h \in L^1(Q_0)_+$. Define
$$h_n = h \cdot \mathbb{1}_{\{f > 1/n\}} \wedge n.$$
Since $f > 0$, $h_n \uparrow h$. Furthermore, $h_n \leq n^2 \cdot f$, because on $\{f > 1/n\}$, $h_n \leq n$. □

Since $h_n \uparrow h$ implies that $S(h_n) \uparrow S(h)$, (A1) can be seen as an approximation property for $M = S(\Delta_{Q_0})$. We will need a similar approximation property for $M^\circ = \{g \in L^1(P_0)_+ : S^*(g) = 1\}$:

(A2) For all $g \in L^1(P_0)_+$ such that $S^*(g) = 1$ and $g > 0$, there exists a sequence $g_n \in L^1(P_0)_+$ such that $S^*(g_n) = \|g_n\|_1 \cdot 1$, $g_n \lesssim g$ and $g_n \uparrow 1$.

We will first show how these two assumptions are used to prove our main theorem, after which we will show in two examples how one checks these assumptions.

THEOREM 3.13. *Let $S : L^1(Q_0) \to L^1(P_0)$ be a coarsening satisfying (A1) and (A2). Then the CAR assumption cannot be tested, so $\mathcal{P}_{\mathrm{CAR}}$ is dense in $\Delta_{P_0}$.*

PROOF. Define $M = S(\Delta_{Q_0})(\subset \Delta_{P_0})$ and $M^\circ = \{g \in L^1(P_0)_+ : S^*(g) = 1\}$. We have to prove that $M \cdot M^\circ$ is dense in $\Delta_{P_0}$. Let $f \in \Delta_{P_0}$ such that $f > 0$ and $f$ is bounded; the set of all these functions is clearly dense in $\Delta_{P_0}$, so it is enough to prove that $f \in \overline{M \cdot M^\circ}$.

Clearly, $\mathrm{KL}_f(1) < +\infty$, so the infimum of $\mathrm{KL}_f$ on $\overline{M}^\sigma \subset E_{P_0}$ is finite (because $1 \in M$). As noted before, since $\mathrm{KL}_f$ is weakly lower semi-continuous and $\overline{M}^\sigma$ is weakly compact, $\mathrm{KL}_f$ attains its minimum somewhere in $\overline{M}^\sigma$, let us say in $k \in \overline{M}^\sigma$. Using Lemma 3.11, we can see that $k_{/\!/} > 0$, since otherwise $\mathrm{KL}_f(k) = \mathrm{KL}_f(k_{/\!/}) = +\infty$ (here we use that $f > 0$). Since $\overline{M}^\sigma = S(\overline{\Delta_{Q_0}}^\sigma)$, we can choose $h_0 \in E_{Q_0,+}$ with $S(h_0) = k$ and $h_{0,/\!/} > 0$ [here we use (A1)].

Let $h \in \Delta_{Q_0}$. According to Lemma 3.12, there exists a sequence $h_n \in L^1(Q_0)_+$ with $h_n \lesssim h_{0,/\!/}$ and $h_n \uparrow h$ [and, therefore, $S(h_n) \uparrow S(h)$]. Define $a_n = \langle h_n, 1 \rangle$ (so $a_n \uparrow 1$) and fix $n$. Because $h_n \lesssim h_{0,//}$, there exists $0 < \varepsilon \leq 1$ such that for all $\lambda \in\, ]-\varepsilon, \varepsilon[$,

$$h_{0,//} + \lambda(a_n^{-1} h_n - h_{0,//}) \geq 0.$$

We conclude that $h_0 + \lambda(a_n^{-1} h_n - h_0) \in \overline{\Delta_{Q_0}}^\sigma$, and so $k + \lambda(a_n^{-1} S(h_n) - k) \in \overline{M}^\sigma$. Therefore, for all $\lambda \in\, ]-\varepsilon, \varepsilon[$,

$$\int -\log\left(\frac{k_{//} + \lambda(a_n^{-1} S(h_n) - k_{//})}{f}\right) f\, dP_0 \geq \int -\log\left(\frac{k_{//}}{f}\right) f\, dP_0,$$

which by differentiating at $\lambda = 0$ implies

$$\left\langle S(h_n), \frac{f}{k_{//}} \right\rangle = a_n.$$

Since $S(h_n) \uparrow S(h)$, we conclude for every $h \in \Delta_{Q_0}$ that $\langle S(h), f/k_{//} \rangle = 1$, which proves that $f/k_{//} \in L^1(P_0)_+$ and that $S^*(f/k_{//}) = 1$, so $f/k_{//} \in M^\circ$.

We would like to conclude that $k_{//} \in \overline{M}$, but we only know that $k \in \overline{M}^\sigma$. We do know, however, that for all $g' \in M^\circ$,

(3.4) $$\langle k_{//}, g' \rangle \leq 1.$$

For if we choose a sequence $g_n \in L^\infty(P_0)_+$ such that $g_n \uparrow g'$, then $S^*(g_n) \leq S^*(g') = 1$, so

$$\langle k_{//}, g' \rangle = \lim \uparrow \langle k_{//}, g_n \rangle \leq \lim \uparrow \langle k, g_n \rangle \leq 1.$$

Now we can use (A2). Define $g = f/k_{//} \in M^\circ$. Clearly, $g > 0$, so there exists a sequence $g_n \in L^1(P_0)_+$ with $S^*(g_n) = \|g_n\|_1 \cdot 1$, such that $g_n \lesssim g$ and $g_n \uparrow 1$. Define $b_n = \|g_n\|_1 \uparrow 1$. There exists $\varepsilon > 0$ such that for all $\lambda \in\, ]-\varepsilon, \varepsilon[$, $g + \lambda(b_n^{-1} g_n - g) \geq 0$, so $g + \lambda(b_n^{-1} g_n - g) \in M^\circ$. Since we have (3.4) and $\langle k_{//}, g \rangle = 1$, we conclude that $\langle k_{//}, g_n \rangle = b_n$, so $\langle k_{//}, 1 \rangle = 1$. But this means that $k = k_{//}$, since $k_{//} \leq k$ and $\langle k, 1 \rangle = 1$. So $k \in L^1(P_0)$ and $k$ is the weak limit [with respect to the duality with $L^\infty(P_0)$] of functions in $M$. However, $L^\infty(P_0)$ is the dual of $L^1(P_0)$ and $M$ is convex, so the weak closure of $M$ in $L^1(P_0)$ equals the strong closure $\overline{M}$, which means that $k \in \overline{M}$. Now choose $\{k_m\} \in M$ such that $\|k - k_m\|_1 \to 0$. We know that $g \cdot k_m \in \Delta_{P_0}$ and $g \cdot k = f \in \Delta_{P_0}$. This means that $\sqrt{gk_m}$ and $\sqrt{gk}$ are positive elements of the unit sphere of $L^2(P_0)$. Since the unit ball is weakly compact in $L^2(P_0)$, we can choose a weakly converging subsequence of $\{\sqrt{gk_m}\}$, let us say $\sqrt{gk_{m_n}} \to \phi$, for some $\phi$ in the unit ball. This means, in particular, that for any $\psi \in L^\infty(P_0)$,

$$\langle \phi, \psi \rangle = \lim_{n \to \infty} \langle \sqrt{gk_{m_n}}, \psi \rangle = \lim_{n \to \infty} \langle \sqrt{k_{m_n}}, \psi\sqrt{g} \rangle = \langle \sqrt{k}, \psi\sqrt{g} \rangle.$$

The last equality follows from the well-known fact that the Hellinger metric induces the same topology on $\Delta_{P_0}$ as the $L^1$-norm, so $\sqrt{k_m} \to \sqrt{k}$ in $L^2(P_0)$, and the fact that $\psi\sqrt{g} \in L^2(P_0)$. Since $L^\infty(P_0)$ is dense in $L^2(P_0)$, we have shown that $\phi = \sqrt{gk}$. This means that every weakly convergent subsequence has the same limit $\sqrt{gk}$, which in turn proves that $\sqrt{gk_m}$ converges weakly to $\sqrt{gk}$. Now note that if for some $\phi \in L^2(P_0)$, we have that $\|\phi\|_2 = 1$, then a neighborhood base for the $L^2$-topology on the unit ball around $\phi$ is given by

$$U_n = \left\{\psi \in L^2(P_0) : \|\psi\|_2 \leq 1 \text{ and } \langle \psi, \phi \rangle > 1 - \frac{1}{n}\right\},$$

since one easily checks that for any $\psi \in U_n$, $\|\psi - \phi\|_2^2 \leq 2/n$. This means that if $\{\psi_n\}$ is a sequence in the unit ball converging weakly to $\phi$, then $\psi_n \to \phi$ in $L^2(P_0)$. Through this we conclude that $\sqrt{gk_m} \to \sqrt{gk}$ in $L^2(P_0)$, which implies that $gk_m \to gk$ in $L^1(P_0)$. So $k \cdot g = f \in \overline{M \cdot M^\circ}$, which proves that $\mathcal{P}_{\text{CAR}}$ is $L_1$-dense in $\Delta_{P_0}$, and, hence, the CAR assumption is not testable. □

We wish to stress that in our opinion the only natural (necessary and sufficient) condition on $S$ for $\mathcal{P}_{\text{CAR}}$ to be dense in $\Delta_{P_0}$ is equation (3.3):

$$S(h) \geq 0 \implies \exists \tilde{h} \geq 0 : S(\tilde{h}) = S(h) \qquad [\forall h \in L^1(Q_0)].$$

This is illustrated by Figure 1 and we have not been able to find counterexamples to this claim. The stronger condition (A1) and condition (A2) were necessary to make our proof work, but must be seen as regularization conditions. We know of examples where (A1) and/or (A2) fail, but we still have the result that CAR cannot be tested. In these examples, the main ideas of the proof of Theorem 3.13 still work, but the details are a bit different.

We will try to illustrate the theorem by two examples, which we will discuss in detail.

EXAMPLE 3.14 (Missing data). Let $Y \in \mathcal{Y}$ be the variable we wish to observe, distributed according to $Q_0$. However, sometimes we can observe $Y$ directly and sometimes the observation is missing, which we will denote by saying that our observation is †. To make things precise, we define our data space $\mathcal{X} = \mathcal{Y} \sqcup \{\dagger\}$. Furthermore, we will use a hidden space to define our coarsening $S$: define $\mathcal{Z} = \mathcal{Y} \times \{0,1\}$ and the map $\psi : \mathcal{Z} \to \mathcal{X}$ as $\psi(y,1) = y$ ($Y$ is not missing), $\psi(y,0) = \dagger$ ($Y$ is missing). Choose $\mu_0 = Q_0 \times (\frac{1}{2}\delta_0 + \frac{1}{2}\delta_1)$, so one possible CAR distribution is that each observation has probability $\frac{1}{2}$ of being missing, independently of $Y \sim Q_0$. It also means that $P_0 = \psi(\mu_0) = \frac{1}{2}\mathbb{1}_\mathcal{Y} \cdot Q_0 + \frac{1}{2}\delta_\dagger$. Then for $h \in L^1(Q_0)$, we define $S(h) \in L^1(P_0)$ as follows:

$$S(h)(y) = E_{\mu_0}(h(Y)|X = y) = h(y) \qquad (\text{for } y \in \mathcal{Y})$$

and

$$S(h)(\dagger) = E_{\mu_0}(h(Y)|X = \dagger) = E_{Q_0}(h(Y)).$$

It is not hard to check that for $g \in L^1(P_0)$,

$$S^*(g)(y) = E_{\mu_0}(g(X)|Y=y) = \tfrac{1}{2}g(y) + \tfrac{1}{2}g(\dagger),$$

so indeed $S^*(1) = 1$, which shows that $S$ is a coarsening. Then

$$\mathcal{P}_{\text{CAR}} = \{g \cdot S(h) : h \in \Delta_{Q_0}, g \geq 0, S^*(g) = 1\}.$$

Since $S^*(g) = 1$ implies that $g(y) = 2 - g(\dagger)$ for $Q_0$-almost all $y$, we see that we get all distributions in $\mathcal{P}_{\text{CAR}}$ by allowing $Y$ to be distributed according to an arbitrary density $h$ with respect to $Q_0$ and assuming that each observation has an arbitrary probability $p = \tfrac{1}{2}g(\dagger)$ to be missing, independently of $Y$.

Now we would like to check assumptions (A1) and (A2). In this case (A2) is trivial, because if $g \in L^1(P_0)$ and $g > 0$ such that $S^*(g) = 1$, then $0 < \min(2 - g(\dagger), g(\dagger)) \leq g$, so $1 \lesssim g$. Assumption (A1) is also not so hard to check, since if we restrict $S(h)$ to $\mathcal{Y}$, we get that $S(h) = \tfrac{1}{2}h$, seen as elements of $L^1(Q_0)$. This shows that $S(h)_{//} = \tfrac{1}{2}h_{//}$, so $S(h)_{//} > 0$ clearly implies $h_{//} > 0$. Theorem 3.13 now states that the CAR assumption cannot be tested in this case, so $\mathcal{P}_{\text{CAR}}$ is dense in $\Delta_{P_0}$. Clearly, in this simple example it is very easy to directly verify that, in fact, $\mathcal{P}_{\text{CAR}} = \Delta_{P_0}$.

EXAMPLE 3.15 (Right-censored data). Let $\mathcal{Y} = ]0,1[$, $Q_0(dt) = dt$ on $\mathcal{Y}$ and $Y$ be a time of interest distributed according to a density with respect to $Q_0$. All that follows can be easily generalized to an arbitrary measure on an open subset of $]0,\infty[$, at the cost of some notational difficulty. Let $C \in ]0,1[$ be a censoring time and let the data $(X, \Delta)$ consist of

$$(X, \Delta) = \psi(Y, C) \stackrel{\text{def}}{=} (Y \wedge C, \mathbb{1}_{\{Y \leq C\}}).$$

We will construct our coarsening $S$ as follows: define $\mu_0 = dt\,dc$ on $\mathcal{Y} \times ]0,1[$. Then define $P_0 = \psi(\mu_0)$ as a probability measure on the data-space $\mathcal{X}$. One easily checks that

$$P_0(dx, \delta) = (1-x)\,dx \cdot \mathbb{1}_{\{\delta=1\}} + (1-x)\,dx \cdot \mathbb{1}_{\{\delta=0\}}.$$

Now define for $h \in L^1(Q_0)$,

$$S(h)(x, \delta) = E_{\mu_0}(h(Y)|(X, \Delta) = (x, \delta)).$$

This is just saying that $S(h)$ is the density of the distribution of the data with respect to $P_0$ when $Y$ and $C$ are independent, $Y$ distributed according to $h(t)\,dt$ and $C$ distributed according to $dc$. Therefore, one easily calculates

$$S(h)(x, \delta) = h(x) \cdot \mathbb{1}_{\{\delta=1\}} + \frac{1}{1-x}\int_x^1 h(t)\,dt \cdot \mathbb{1}_{\{\delta=0\}}.$$

Clearly, $S(1) = 1$ and $S$ is positive. Furthermore, for all $h \in L^1(Q_0)$ and $g \in L^\infty(P_0)$,

$$\langle S(h), g \rangle_{P_0} = \int_0^1 (1-x)h(x)g(x,1)\,dx + \int_0^1 \int_x^1 h(t)g(x,0)\,dt\,dx$$
$$= \int_0^1 (1-t)g(t,1)h(t)\,dt + \int_0^1 \left( \int_0^t g(x,0)\,dx \right) h(t)\,dt$$
$$= \int_0^1 \left( (1-t)g(t,1) + \int_0^t g(x,0)\,dx \right) h(t)\,dt,$$

so we see that

$$S^*(g)(t) = (1-t)g(t,1) + \int_0^t g(x,0)\,dx.$$

Hence, $S^*(1) = 1$, so $S$ is indeed a coarsening. Define $M = S(\Delta_{Q_0})$ and $M^\circ = \{g \in L^1(P_0)_+ : S^*(g) = 1\}$. Let $h \in \Delta_{Q_0}$ and $g \in M^\circ$. Since $S^*(g) = 1$,

$$g(t,1) = \frac{1}{1-t}\left( 1 - \int_0^t g(x,0)\,dx \right).$$

Because $g \geq 0$, we have that $\int_0^1 g(x,0)\,dx \leq 1$. If $\int_0^1 g(x,0)\,dx = 1$ and we let $C$ be distributed according to $g(x,0)$ and $Y$ according to $h$, we can easily check that the density of $\psi(Y,C)$ with respect to $P_0$ is exactly $S(h) \cdot g$, so $\mathcal{P}_{\mathrm{CAR}}$ contains all data distributions one gets if $Y$ and $C$ are independent and dominated by the Lebesgue measure. If we allow $C$ to be distributed according to a subdensity (i.e., just saying that the censoring time has a positive probability of being bigger than the largest possible value for $Y$), then we get all of $\mathcal{P}_{\mathrm{CAR}}$.

We would now like to check assumptions (A1) and (A2). Denote by $S(h)_{\{\delta=1\}}$ the restriction of $S(h)$ to $\{\delta = 1\}$. Clearly, for $h \in E_{Q_0}$, $S(h)_{\{\delta=1\}} = (1-y) \cdot h$; here $(1-y) \cdot h$ acts on $\phi(y) \in L^\infty(Q_0)$ as follows:

$$\langle (1-y) \cdot h, \phi(y) \rangle = \langle h, (1-y) \cdot \phi(y) \rangle.$$

To check (A1) it is enough to conclude that $h_{//} > 0$ whenever $[(1-y) \cdot h]_{//} > 0$ and $h \geq 0$. However, in that case $(1-y) \cdot h \leq h$ (because $\|1-y\|_\infty \leq 1$), so $0 < [(1-y) \cdot h]_{//} \leq h_{//}$.

Now let $g \in M^\circ$, $g > 0$. Define

$$g_n(t,0) = \mathbb{1}_{\{g(t,0) > 1/n\}}.$$

Clearly, $g_n(t,0) \leq ng(t,0)$. Define $\lambda_n = \int_0^1 g_n(x,0)\,dx$ and

$$g_n(t,1) = \frac{1}{1-t}\left( \lambda_n - \int_0^t g_n(x,0)\,dx \right).$$

Then we have that

$$g_n(t,1) = \frac{1}{1-t}\int_t^1 g_n(x,0)\,dx$$
$$\leq \frac{n}{1-t}\int_t^1 g(x,0)\,dx$$
$$= \frac{n}{1-t}\left(\int_0^1 g(x,0)\,dx - \int_0^t g(x,0)\,dx\right)$$
$$\leq \frac{n}{1-t}\left(1 - \int_0^t g(x,0)\,dx\right)$$
$$= ng(t,1).$$

So $g_n \lesssim g$. By construction we have that $S^*(g_n) = \lambda_n \cdot 1$. Furthermore, since $g > 0$, $g_n(t,0) \uparrow 1$. This implies that $g_n(t,1) \uparrow 1$, so $g_n \uparrow 1$. This proves that assumption (A2) is also satisfied. Theorem 3.13 now states that the CAR assumption cannot be tested in the case of right-censored data. We wish to remark that this in itself is not a new result, but merely an illustration of Theorem 3.13.

## APPENDIX

In this Appendix we will give the proof of a lemma which is a bit technical. We repeat some notation: define for a probability measure $P$ the space $E = (L^\infty(P))'$. This is an ordered vector space and $L^1(P)$ is a band in $E$, which means that each $h \in E_+$ can be uniquely decomposed as $h = h_{//} + h_\perp$, where $h_{//} \in L^1(P)_+$ and $h_\perp \geq 0$ is disjoint from $P$, so for each $f \in L^1(P)_+$, $\inf(f, h_\perp) = 0$. According to Schaefer and Wolff [(1999), Chapter V, Theorem 1.5] this is equivalent to saying that for each $\phi \in L^\infty(P)_+$ and each $\varepsilon > 0$, there exists a decomposition $\phi = \phi_1 + \phi_2$, $\phi_1 \geq 0$, $\phi_2 \geq 0$, such that $\langle h_\perp, \phi_1 \rangle + \langle f, \phi_2 \rangle < \varepsilon$. For convenience, we repeat the definition of $\mathrm{KL}_f$ for $f \in L^1(P)_+$. Define for $h \in E_+$,

$$\mathrm{KL}_f(h) = \sup\left\{\sum_{i=1}^n -\log\left(\frac{\langle h, \phi_i \rangle}{\langle f, \phi_i \rangle}\right)\langle f, \phi_i \rangle : \phi_i \in L^\infty(P)_+, \sum_{i=1}^n \phi_i = 1\right\}.$$

LEMMA A.1. *Let $f \in L^1(P)_+$. Then, in the notation introduced above, for all $h \in E_+$,*

$$\mathrm{KL}_f(h) = \mathrm{KL}_f(h_{//}) = \int -\log\left(\frac{h_{//}}{f}\right) f\,dP.$$

PROOF. The second equality is well known for the Kullback–Leibler divergence [see, e.g., Pinsker (1964), Section 2.4] and can be proved using standard techniques like monotone classes. As for the first, since $h \geq h_{//}$, it is clear that $\mathrm{KL}_f(h_{//}) \geq \mathrm{KL}_f(h)$ ($-\log$ is a decreasing function). Assume $\mathrm{KL}_f(h) < +\infty$. We would like to make an important observation: if $\sum_{i=1}^n \phi_i = 1$, and we decompose each $\phi_i = \phi_{i,1} + \phi_{i,2}$, then Jensen gives us

$$\text{(A.1)} \quad \sum_{i=1}^n \sum_{j=1}^2 -\log\left(\frac{\langle h, \phi_{i,j}\rangle}{\langle f, \phi_{i,j}\rangle}\right)\langle f, \phi_{i,j}\rangle \geq \sum_{i=1}^n -\log\left(\frac{\langle h, \phi_i\rangle}{\langle f, \phi_i\rangle}\right)\langle f, \phi_i\rangle.$$

We only need to consider $\phi_i$ such that $\langle f, \phi_i\rangle > 0$ [we define $\log(1/0) \cdot 0 = 0$]. We also know that for each such $\phi_i$, $\langle h_{//}, \phi_i\rangle > 0$. For if not, we could decompose $\phi_i = \phi_{i,1} + \phi_{i,2}$, such that $\langle h, \phi_{i,1}\rangle$ is arbitrarily small and $\langle f, \phi_{i,1}\rangle > \langle f, \phi_i\rangle/2$, which by (A.1) would imply that $\mathrm{KL}_f(h) = +\infty$.

So consider $\phi_i \geq 0$ with $\langle h_{//}, \phi_i\rangle > 0$ and $\mathbb{1}_{\{f>0\}} \leq \sum_{i=1}^n \phi_i \leq 1$. Let $\varepsilon > 0$. Because $\inf(h_\perp, f) = \inf(h_\perp, h_{//}) = 0$, we can find a decomposition $\phi_i = \phi_{i,1} + \phi_{i,2}$ for each $i$ such that $\langle h_{//}, \phi_i - \phi_{i,1}\rangle = \langle h_{//}, \phi_{i,2}\rangle < \delta$, $\langle f, \phi_i - \phi_{i,1}\rangle = \langle f, \phi_{i,2}\rangle < \delta$ and $\langle h_\perp, \phi_{i,1}\rangle < \delta$. Here we can choose $\delta > 0$ such that

$$\sum_{i=1}^n -\log\left(\frac{\langle h_{//} + h_\perp, \phi_{i,1}\rangle}{\langle f, \phi_{i,1}\rangle}\right)\langle f, \phi_{i,1}\rangle \geq \sum_{i=1}^n -\log\left(\frac{\langle h_{//}, \phi_i\rangle}{\langle f, \phi_i\rangle}\right)\langle f, \phi_i\rangle - \varepsilon$$

and, noting that $\langle h, 1\rangle > 0$, since $\mathrm{KL}_f(h) < +\infty$,

$$\sum_{i=1}^n -\log\left(\frac{\langle h, 1\rangle}{\langle f, \phi_{i,2}\rangle}\right)\langle f, \phi_{i,2}\rangle \geq -\varepsilon.$$

This last inequality implies that

$$\sum_{i=1}^n -\log\left(\frac{\langle h, \phi_{i,2}\rangle}{\langle f, \phi_{i,2}\rangle}\right)\langle f, \phi_{i,2}\rangle \geq \sum_{i=1}^n -\log\left(\frac{\langle h, 1\rangle}{\langle f, \phi_{i,2}\rangle}\right)\langle f, \phi_{i,2}\rangle \geq -\varepsilon.$$

All in all, we can conclude that

$$\mathrm{KL}_f(h) \geq \sum_{i=1}^n \sum_{j=1}^2 -\log\left(\frac{\langle h, \phi_{i,j}\rangle}{\langle f, \phi_{i,j}\rangle}\right)\langle f, \phi_{i,j}\rangle$$

$$\geq \sum_{i=1}^n -\log\left(\frac{\langle h_{//}, \phi_i\rangle}{\langle f, \phi_i\rangle}\right)\langle f, \phi_i\rangle - 2\varepsilon.$$

This proves that $\mathrm{KL}_f(h) \geq \mathrm{KL}_f(h_{//})$. □

**Acknowledgments.** The ideas for this paper were mainly conceived while I was a guest at the Institut Henri Poincaré to visit the Statistical Semester. I would like to thank the members of this institute and, in particular, Professor Lucien Birgé for giving me this opportunity.

DEPARTMENT OF MATHEMATICS
FACULTY OF INFORMATION TECHNOLOGY
  AND SYSTEMS
DELFT UNIVERSITY OF TECHNOLOGY
P.O. BOX 5031
2600 GA, DELFT
THE NETHERLANDS
E-MAIL: e.a.cator@ewi.tudelft.nl